\documentclass[preprint,12pt]{elsarticle}

\usepackage{amssymb}

\journal{arXiv}

\newtheorem{thm}{Theorem}
\newtheorem{lem}{Lemma}
\newtheorem{rem}{Remark}
\newtheorem{df}{Definition}
\newtheorem{pr}{Proposition}
 \newproof{pf}{Proof}

\begin{document}

\begin{frontmatter}



\title{ T1 theorem for Campanato spaces on  domains}


\author{Andrei V.\ Vasin\fnref{label2}}

\fntext[label2]{The author was supported by RFBR [grant number 17-01-00607a]}

\ead{andrejvasin@gmail.com}

\address{Admiral Makarov State University of Maritime and Inland Shipping, 198255,
Dwinskaya Street. 5/7, St.-Petersburg, Russia}

\begin{abstract}

Given a   Lipschitz domain $D\subset \mathbb{R}^d,$   a Calder\'{o}n-Zygmund operator $T$ and a modulus of continuity $\omega(x),$ we solve a problem when   the restricted operator $T_Df=T(f\chi_D)\chi_D$ sends the   Campanato space $\mathcal{C}_\omega(D)$ into itself.

 The solution is a T1 type sufficient and necessary condition for the characteristic function $\chi_D$ of $D$:
    $$(T\chi_D)\chi_D \in \mathcal{C}_{\tilde{\omega}}(D),$$
assumed $\tilde{\omega}(x)= \omega(x)/\int_x^1 \omega(t)dt/t.$

To check the hypotheses of T1 theorem we need extra restrictions  on both the boundary of $D$ and   the operator $T.$ It is proved that the restricted Calder\'{o}n-Zygmund operator $T_D$ with the even kernel is bounded on $\mathcal{C}_\omega(D),$ provided $D$ be  $C^{1,\tilde{\omega}}-$smooth domain. This result is sharp.

\end{abstract}


\end{frontmatter}


\section{Introduction.}

\subsection{Restricted Calder\'{o}n-Zygmund operators.} A smooth  homogeneous Calder\'{o}n-Zygmund operator is a principal
value convolution operator of type
$$Tf(y)=PV\int f(x) K(y-x) dx,$$
 where $dx$ denotes the Lebesgue measure in $\mathbb{R}^d$ and
$$ K(x) =\frac{\Omega(x)}{|x|^d },\quad x \neq 0,$$
$\Omega(x)$ being a homogeneous function of degree 0, continuously differentiable  on
$\mathbb{R}^d \setminus 0$ and with zero integral on the unit sphere. Given a  domain $D \subset\mathbb{R}^d$ we consider a modification of $T.$ An  operator  defined by the formula
$$T_Df=  T(\chi_Df)\chi_D$$
is called to be a \textit{restricted} Calder\'{o}n-Zygmund operator.

 \subsection{Moduli of continuity and Campanato spaces.} A continuous   increasing  function $\omega:[0,\infty)\rightarrow[0,\infty),\quad \omega(0)=0$ or $\omega\equiv 1$
   is called a modulus of continuity. In what follows we consider moduli of continuity,  such that a function
 $\frac{\omega(x)}{x^{\varepsilon}}$ is almost decreasing for some  $\varepsilon<1,$  that is
 $$\frac{\omega(s)}{s^{\varepsilon}} <C\frac{\omega(t)}{t^{\varepsilon}}, \quad 0<t<s$$
 with universal constant $C>0$ ( in \cite{Ja2} this one is called  a growth function of upper type less then 1)

  Let $Q$ be a cube in $\mathbb{R}^d$ with edges parallel to coordinate axes. In what follows we denote by  $|Q|$ the volume of $Q$ and by $\ell=\ell(Q)$   its side length.
   Given a modulus of continuity $\omega$  the homogeneous Campanato space $\mathcal{C}_{\omega}(D)$  in arbitrary domain $D\subset  \mathbb{R}^d$ is defined by the seminorm

\begin{equation}\label{eq:eq1}
   \|f\|_{\omega,D}=\sup_{Q\subset D} \inf_{b_Q\in\mathbb{C}}\frac{1}{\omega(\ell)} \|f-b_Q\|_{L^1(Q, dx/|Q|)}.
\end{equation}
 Here  we have $\mathcal{C}_{\omega}(D)= BMO(D)$ when $\omega\equiv 1.$  In the case  $\omega(t)=t^\alpha, \; 0<\alpha<1$ one has $\mathcal{C}_{t^\alpha}(D)=Lip_\alpha (D)$ (Campanato\cite{C}, Meyer\cite{M}).The arguments \cite{C,M} (see, also \cite{KK}, Ch.4) used in the studies of the
classical space $\textrm{BMO}(\mathbb{R}^d)$,  guarantee that, we may replace the $L^1$- norm in (\ref{eq:eq1}) by $L^p$-norm with $1\leq p<\infty$ in arbitrary domain $D.$ This fact will be used in Section 3.   For a bounded Lipschitz domain $D$ we prove below in Section 2 that  $\mathcal{C}_{\omega}(D)\subset L^1(D,dx).$ This enables us to provide the space  $\mathcal{C}_{\omega}(D)$ with a norm $\|f\| =\|f\|_{\omega,D}+\|f\|_{L^1(D,dx)}.$

The spaces  $\mathcal{C}_{\omega}(\mathbb{R}^d)$  are known to be invariant under certain smooth convolution
 Calder\'{o}n-Zygmund operators (Peetre \cite{P}, Jansson \cite{Ja2,Ja1}). In the setting of function spaces defined on domains $D\subset\mathbb{R}^d,$ the following result of Mateu, Orobitg and Verdera ( Main Lemma \cite{MOV}) is crucial.
 \begin{thm}
\label{thm1}

   Let $D$ be a   bounded domain with the  $C^{1,\alpha}-$ smooth boundary, $0<\alpha<1.$ Then the smooth  homogeneous Calder\'{o}n-Zygmund operator $T$ with the even kernel sends $Lip_\alpha(D)$ in $Lip_\alpha(D).$

\end{thm}

The proof   contains  a certain T1 theorem as a essential part, when one checks the boundedness of operator on the characteristic function of $D.$

 The  purpose of the paper it to extend the above assertion on the gap between $Lip_\alpha(D)$ and $BMO(D).$ We do this in the scale of the  Campanato spaces. Before formulating the main result, define  another  modulus associated  with the previous one,
 \begin{equation}\label{eq:eq2}
   \tilde{\omega}(x)= \frac{\omega(x)}{\int_x^1 \omega(t)dt/t}.
    \end{equation}
 It is obvious that $\tilde{\omega}$ is more smooth then $\omega.$ If a modulus of continuity $\omega$ is Dini regular,  that is the  integral
 $$\int_0 \frac{\omega(t)}{t}dt$$
converges, then $\tilde{\omega}(x)\approx \omega(x).$

  We prove the following non symmetric T1 theorem:
 \begin{thm}
\label{thm2}
   Let $\omega$ be a    modulus of continuity (or $\omega\equiv 1$), $D\subset \mathbb{R}^d$ be a bounded  Lipshitz domain and let $\chi_D$ be a characteristic function of the domain $D.$ Let $T$ be  a smooth  homogeneous Calder\'{o}n-Zygmund operator.

    Then the \textit{restricted} operator $T_D$ is bounded on the Campanato space $\mathcal{C}_{\omega} (D),$ if and only if:
   \begin{equation}\label{eq:eq3}
   (T\chi_D)\chi_D \in \mathcal{C}_{\tilde{\omega}}(D).
    \end{equation}
 \end{thm}

Several remarks are in order.
\begin{rem}
 \label{rem1}
 Particularly, we have that (\ref{eq:eq3}) coincides with
 $$(T\chi_D)\chi_D \in \mathcal{C}_{\omega}(D),$$ provided $\omega$ is Dini regular.
\end{rem}
 \begin{rem}
 \label{rem2}
  If $\omega_{\alpha}(t)=\log^{-\alpha}t^{-1}, \: 0\leq \alpha<1,$ then  $\tilde{\omega}_{\alpha}(t)\approx \log^{-1}t^{-1},$ and therefore we have the common necessary and sufficient condition for all moduli $\omega_{\alpha}$
 $$(T\chi_D)\chi_D \in \mathcal{C}_{\log^{-1}t^{-1}}(D).$$
 \end{rem}
 \begin{rem}
 \label{rem3}
 For any cube $Q$ one has that   $T\chi_Q\in BMO(\mathbb{R}^d),$  but $BMO(Q)$ is not invariant under action of $T_Q.$
\end{rem}

We assume  that the underlaying domain is Lipschitz and do not consider more general domains.
 In fact, all that we need  is to have possibility to extend functions from the   Campanato spaces defined  on domains  to ambient space $\mathbb{R}^d.$  Further in order to verify  (\ref{eq:eq3}),   we claim extra smoothness of the boundary of the underlaying domain then to be only Lipschitz. Therefore, now we state on the Lipschitz domain and no more.

  An extension of functions from the Campanato space reminds us the Stein   (\cite{S} Ch. 6) about extension of the functions from  Lipshitz spaces $\textrm{Lip}_{\omega}(D).$

 In the present paper  we apply  approach  on extension in BMO (Jones \cite{Jo}) and Besov spaces (DeVore, Sharpley \cite{DVS}) on domains. To do this job, we  define a variant of Campanato
 space $\mathcal{C}^{int}_{\omega}(D)$  distinct from the previous one. We take supremum  in  the  seminorm below with respect to
 all cubes  $Q\subset D,$    separated from the boundary $\partial D:$
\begin{equation}\label{eq:eq4}
   \|f\|^{int}_{\omega,D}=\sup_{Q: 2Q\subset D} \inf_{b_Q\in\mathbb{C}}\frac{1}{\omega(\ell)} \|f-b_Q\|_{L^1(Q\cap D, dx/|Q|)}.
\end{equation}
 In fact, for the Lipshitz domain $D$ the seminorms  (\ref{eq:eq1}) and  (\ref{eq:eq4}) are equivalent and define the same space $\mathcal{C}^{int}_{\omega}(D)=\mathcal{C}_{\omega}(D).$ Moreover,  it is just this case to extend the functions from the space $\mathcal{C}^{int}_{\omega}(D)$ into the ambient space $\mathbb{R}^d.$
The precise proposition see in Extension Lemma \ref{lem1}  below in Section 2.

\subsection{Applications of T1 theorem.}

  In this section we restrict our attention to operators with an even kernel, that is $K(-x)=K(x).$ Also we assume that underlaying domain has  the $C^{1,\omega}$-smooth  boundary (Lipschitz domains with $C^{1,\omega}$-smooth parametrization of the boundary). Recall, that
$\omega\equiv 1$ is assumed to be  the modulus of continuity. In a  following sharp result we establish  the relationship  of the smoothness properties of the boundary of a domain  with the boundedness of the restricted Calder\'{o}n-Zygmund operator in  the   Campanato spaces.
\begin{thm}
\label{thm3}
Let $\omega$ be a   modulus of continuity, $\tilde{\omega}$ is given by (\ref{eq:eq2}), $D\subset \mathbb{R}^d$ be a  bounded domain with the $C^{1,\tilde{\omega}}$-smooth  boundary.  Let $T$ be  an even  smooth  homogeneous Calder\'{o}n-Zygmund operator.  Then  $T_D$ is bounded on the Campanato space  $\mathcal{C}_{\omega} (D).$

\end{thm}
Theorem \ref{thm3} is easy consequence of a next result.
\begin{pr}
\label{pr1}
 Given a modulus of continuity $\omega,$   a  domain $D\subset \mathbb{R}^d$ with the $C^{1,\omega}-$smooth  boundary,  let $T$ be  a smooth  homogeneous Calder\'{o}n-Zygmund operator with an even kernel.  Then  $(T\chi_D)\chi_D \in\mathcal{C}_{\omega} (D).$
\end{pr}
We  reduce Proposition \ref{pr1} to a assertion about  the weighted Bloch spaces.
\begin{df}
 \label{df11}
 Given a modulus of continuity $\omega$, and domain $D\subset \mathbb{R}^d,$ the real weighted Bloch space $\mathcal{B}_\omega (D)$ consists of $C^1-$smooth in $D$ functions $f$ such that
 \begin{equation}
 \label{eq:eq41}
 \sup_{x\in D} \frac{|\nabla f(x)\rho(x)|}{\omega(\rho(x))}< \infty,
\end{equation}
with $\rho(x)= dist( x,\partial D).$
\end{df}
In Section 4 we will prove auxiliary Lemma \ref{lem7} about imbedding  $\mathcal{B}_\omega (D)\subset \mathcal{C}_\omega (D).$
This imbedding is sharp. We show this for harmonic functions on  $C^1-$smooth domains  in Lemma \ref{lem8}.   With these results in hand   Theorem \ref{thm3} follows from  the last reduction.

\begin{pr}
\label{pr2}
   Let $\omega$ be a  modulus of continuity   and let $D\subset \mathbb{R}^d$ be a  $C^{1,\omega}-$smooth domain. Let $T$ be  an even  smooth  homogeneous Calder\'{o}n-Zygmund operator. Then   $(T\chi_D)\chi_D\in \mathcal{B}_{\omega} (D).$
\end{pr}
It is worth mentioning that Proposition \ref{pr2} is  also sharp. It is  proved in  \cite{V}  in a particular, but  important case of the Beurling transform in plane domains. The same argument may be modified for arbitrary  homogeneous Calder\'{o}n-Zygmund operator with the smooth even kernel. This remark shows that Theorem \ref{thm3} is sharp too.

\subsection {Organization and notation}
In Section 2 we shall introduce some basic notation and set up some necessary
preliminaries about Whitney covering in order to construct Extension lemma.   Theorem \ref{thm2} is proved in Section 3. Imbedding Lemma and the proof of Theorem \ref{thm3} are postponed  to Section 4.

As usual, the letter $C$ will denote a constant, which may be different at each
occurrence and which is independent of the relevant variables under consideration.
The notation $A\lesssim B$ means that there is a fixed positive constant $C$ such that $A<CB.$ If  $A\lesssim B\lesssim A,$ then we write $A\approx B.$ We write $A\lesssim_{a,b} B,$ if a corresponding constant depends on $a,b.$ We denote by

\section {Extension lemma}
We start with a definition of a Lipschitz domain. We also introduce a notation $sQ$ to denote a cube with the same centre as $Q$ and with side length $s\ell(Q).$
\begin{df}
 \label{df2}
 A bounded domain $D\subset  \mathbb{R}^d$ is called to be a $(\delta, R)-$Lipshitz domain if for  each  $a\in\partial D$ there exist the function $A:\mathbb{R}^{d-1}\rightarrow \mathbb{R}$ with $\|\nabla A\|_\infty\leq \delta$  and the cube $Q$ with side length $R$ and centre $a$ such that, after a suitable shift and rotation
 $$D\cap Q=\{ (x,y)\in (\mathbb{R}^{d-1}, \mathbb{R})\cap Q:y>A(x)\}.$$
This cube is called a $R$-window of domain. Without risk of confusion we will forget about parameters $\delta$ and $R$, and will talk about Lipshitz domains in general.
\end{df}

In what follows we need to introduce a Whitney covering. Consider a given dyadic grid of semi-open cubes in $\mathbb{R}^d$.

 \begin{df}
\label{df3}
We say that a collection of cubes $\mathcal{W}$ is a Whitney covering of a Lipschitz domain $D$ if

\begin{enumerate}
  \item The cubes in $\mathcal{W}$  are dyadic.
  \item The cubes have pairwise disjoint interiors.
  \item The union of the cubes in $\mathcal{W}$ is $D.$
  \item There is an estimate
$$\textrm{diam}(Q)\leq \textrm{dist}(Q, \partial D) \leq 4\textrm{diam}(Q).$$
  \item  Two neighbor cubes $Q$ and $R$ (i.e. $\bar{Q}\cup \bar{R}\neq \emptyset$) satisfy $\ell(q)\leq 2\ell(R)$
  \item  The family of $\{10Q\}_{Q\in\mathcal{W}}$ have finite superposition (i.e. $\sum_{Q\in\mathcal{W}}\chi_{10Q}<\infty$)
\end{enumerate}

\end{df}
We do not prove here the existence of such a covering because this kind of coverings are well
known and widely used in the literature (e.g. Stein, Ch.6\cite{S}). Our definition is taken from  \cite{PT} (Pratz, Tolsa).

Recall that we consider the $R$-window $Q$ to be a cube centered in $x\in \partial D,$
 with side-length $R$
inducing a Lipschitz parameterization of the boundary. Thus, each $Q$ induces a vertical direction,
given by the eventually rotated $x_d$ axis. The following
is an easy consequence of the previous statements and the fact that the domain is Lipschitz:

 7.\emph{The number of Whitney cubes in a window with the same side-length intersecting
a given vertical line is uniformly bounded where the vertical direction is the one induced by the window.}

This is the last property of the Whitney cubes we want to point out. In fact,  we consider   Whitney covering $\mathcal{W}$ of a Lipschitz domain $D$ as well as a  Whitney covering $\mathcal{W}'$ of a complement $D'=\mathbb{R}^d\setminus \textrm{clos}(D).$

 To construct an extension we first define a $C^{\infty}$-smooth partition of unity $\{\psi_Q\}_{Q\in \mathcal{W}'}$  associated with the Whitney covering   $\mathcal{W}'$ of $D'.$ This means that each bump function satisfies following conditions: $\psi_Q$ is $C^{\infty}$-smooth, $\chi_{4/5 Q}\leq \psi_Q\leq \chi_{5/4 Q},\:Q\in \mathcal{W}'$ and $\sum_{Q\in \mathcal{W}'} \psi_Q =\chi_{D'}.$

Given a Whitney cube $Q\in \mathcal{W}',$   define  a  Whitney cube $\tilde{Q}\in \mathcal{W},$  called reflective to $Q$, if it is a maximal cube such that $\textrm{dist}(Q, \tilde{Q})\leq 2 \textrm{dist}(Q, \partial D).$ Let $f_Q=\frac{1}{|Q|}\int_Qfdx$ be a mean value of $f$ over a cube $Q.$

Define an extension  of $f$ by formula:
\begin{equation}\label{eq:eq5}
 \tilde{f}=f\chi_D+\sum_{Q\in \mathcal{W}',\:\ell(Q)\leq R} \psi_Q f_{\tilde{Q}}.
\end{equation}
We are ready to formulate Extension lemma.
\begin{lem}
\label{lem1}
   Given $ f \in \mathcal{C}^{int}_{\omega}(D),$ then the function  $\tilde{f},$ defined by (\ref{eq:eq5}),  satisfies following properties:
   \begin{enumerate}
   \item $\tilde{f}$ has compact support,
    \item $\tilde{f}$ is $C^\infty-$smooth in $D'=\mathbb{R}^d\backslash \textrm{clos}(D),$
    \item $\tilde{f} \in L^1(\mathbb{R}^d,dx)\cap\mathcal{C}_{\omega}(\mathbb{R}^d),$ and  there is an  estimate
     $$\|\tilde{f}\|_{\omega,\mathbb{R}^d}+\|\tilde{f}\|_{L^1(\mathbb{R}^d,dx)}\lesssim\|f\|^{int}_ {\omega,D}+\|f\|_{L^1(D,dx)}$$
   \end{enumerate}
   \end{lem}
  \begin{pf}
  The properties 1) and 2) are obvious. To prove 3) we have to estimate supremum in (\ref{eq:eq1}), provided $D=\mathbb{R}^d.$

  First we consider only small cubes $Q$ in $\mathbb{R}^d$ such that $2Q\cap \partial D\neq \emptyset.$ To define  exactly the value of smallness,  we choose a constant $r_0< R,$ so that the following property, that will be needed along the proof, is satisfied. For each cube $Q\in \mathcal{W}$ let $Q'=9/8Q,$ then $2Q'\in D.$

 \begin{lem}
\label{lem2}
    Let $D$ be a Lipschitz domain. There exist three positive constants $C, \:c,\: r_0, $ depending only on the Lipschitz constants of $D,$ such that if  $Q$ is a cube in $\mathbb{R}^d$ with $\ell(Q)< r_0$ and such that $2Q\cap \partial D\neq \emptyset,$ then
   $$\int_Q|\tilde{f}-\tilde{f}_Q|dx\leq C \sum_{S\subset cQ,\:S\in \mathcal{W}}\int_{S'}|f-f_{S'}|dx.$$

\end{lem}
 This is mostly Lemma 5.2 \cite{DVS}, which   is proved  for a wider class of  domains. It is formulated in \cite{DVS} in a position  that $f\in L^1(D, dx).$ The same argument
 takes place, when $f$ is summing only over all cubes $ S\in \mathcal{W}.$
 We choose the parameter $r_0,$ such that the cube $cQ$ lies in a certain  R-window  and conditions of Definition \ref{df2}  are satisfied.

 Since for any Whitney cube $S\in \mathcal{W}$ we have that $ 2S'\subset D,$ hence
 \begin{equation}\label{eq:eq6}
 I=\int_Q|\tilde{f}-\tilde{f}_Q|dx\leq C\|f\|_{\omega,D}^{int} \sum_{S\subset cQ,\:S\in \mathcal{W}}\omega (\ell(S)) \ell(S)^d .
\end{equation}
 Now, we estimate the number of cubes of the same size in the sum. The property 7 from Definition \ref{df3} of Whitney covering
  give us that for any Whitney cube $S\subset cQ$ there is a vertical line (defined by axis $x_d$) intersected with
  finitely many of Whitney cubes of the size $\ell(S)$ . This number is estimated from above by a constant
  $C$ depending on Lipshitz constant of domain $D.$ Thus, we may estimate the number of all cubes  of the  size
  $\ell(S)$ intersected with $cQ$ as
\begin{equation}\label{eq:eq7}
   \sharp(S)\lesssim \left(\frac{\ell(cQ)}{\ell(S)}\right)^{d-1}
\end{equation}
 Let $s$ be an integer such that $2^s= \ell(S)$ and let $n$ be an   integer such that $2^n\leq \ell(cQ)<2^{n+1}.$ Then (\ref{eq:eq7}) is estimated by
 $$\sharp(S)\lesssim\left(\frac{2^n}{2^s}\right)^{d-1}$$
 with constant not depending on $Q.$ So, continue one has
 $$I\lesssim\sum_{s=-\infty}^n\left(\frac{2^n}{2^s}\right)^{d-1}\omega (2^s) (2^s)^d=$$
 $$=(2^n)^{d-1}\sum_{s=-\infty}^n{2^s}\omega (2^s) \lesssim$$
 $$\lesssim(2^n)^{d}\omega (2^n) \approx |Q| \omega (\ell(Q)), $$
 and we are done for small cubes near the boundary.

To deal with large cubes $Q,\quad\ell(Q)> r_0,$ we prove that $\tilde{f}\in L^1(\mathbb{R}^d, dx).$
This gives us the claimed estimate for large cubes:

$$I=\frac{1}{|Q|}\int_{Q}|\tilde{f}|dx\leq\frac{1}{r_0^{d}}\|\tilde{f}\|_{ L^1(\mathbb{R}^d, dx)}\lesssim\omega(\ell(Q)).$$
To prove that $\tilde{f}\in L^1(\mathbb{R}^d, dx),$ we first choose a family    $E=\{Q_k\}_{k=1}^{N_0}$ of cubes $Q_k$ with $\ell(Q_k)\approx r_0,$  such that
$$\partial D\subset \cup_{k=1}^{N_0} Q_k.$$
Then $N_0 \thickapprox \mathcal{H}^{d-1}( \partial D)/r_0^{d-1}.$

Second consider a  family $F=\{Q_k\}$ of Whitney cubes from $\mathcal{W}\cup \mathcal{W}'$ that covering  $\textrm{supp}(\tilde{f}) \setminus E.$ It is obvious that $\ell(Q_k)\gtrsim r_0,$  hence the family $F$ is finite.
Now, the function $\tilde{f}$ is summing in each cube  $Q\in E\cup F$ in a finite family, and therefore $\tilde{f}\in L^1(\mathbb{R}^d, dx).$

 Finally, consider small cubes in $D'$ far from the boundary. Since, the function $\tilde{f}$ is $C^\infty$-smooth in $D',$ and the set  $\textrm{supp}\tilde{f}\backslash (D\cup E )$ is compact, hence the claimed estimates (\ref{eq:eq1}) easily follow.

 The rest of the proof is done by the closed graph theorem.$\Box$
   \end{pf}

\section {Proof of T1 theorem}

Recall that  $\omega$ is a    modulus of continuity, $D\subset \mathbb{R}^d$ is a  Lipshitz domain and let $\chi_D$ be a characteristic function of the domain $D.$ Let $T_D$ be  a restricted Calder\'{o}n-Zygmund operator. Fix $f\in\mathcal{C}_{\omega}(D)$ and
 let us  check when $T_Df\in\mathcal{C}_{\omega}(D).$  Given a cube  $Q$  with size length $\ell,$ we have to estimate oscillation
\begin{equation}\label{eq:eq8}
   I=\frac{1}{|Q|}\int_{Q\subset D}|Tf-b_Q|dx < C \omega (\ell )\|f\|
\end{equation}
with  an appropriate constant $b_Q$ and $\|f\| =\|f\|_{\omega,D}+\|f\|_{L^1(D,dx)}.$ By Extension Lemma \ref{lem1}  we may consider only cubes $Q$ such that $2Q\subset D.$ We start with a modification of a construction from the proof of Theorem 1.5 in  \cite{DV}(see, also \cite{H}). Taking into account a mean value $f_Q$ of $f$ over a cube $Q,$ we put
$$f_1=f_Q \chi_D,$$
$$f_2=(f-f_Q) \chi_{2Q},$$
$$f_3=(f-f_Q) \chi_{D\backslash 2Q}.$$
Observe that $f = f_1 + f_2 + f_3$.

\begin{lem}
\label{lem3}
     There exist constant $C>0$ depending on $f$ such that for the  functions $f_k$ with $k=2,3$   one has
   $$\frac{1}{|Q|}\int_{Q}|Tf_k-b_{k,Q}|dx\leq C\omega(\ell),$$
   with an appropriate constants $ b_{k,Q}.$
   \end{lem}
\begin{pf} We need a simple standard auxiliary argument
\begin{lem}
\label{lem4}
    Let $f\in \mathcal{C}_{\omega}(D),$ then

\begin{equation}\label{eq:eq9}
  |f_Q-f_{2Q}|\lesssim \omega(\ell).
\end{equation}
   \end{lem}

a) For $k=2$  we proceed by choosing $b_{2,Q}=0$ and  H\"{o}lder inequality
$$I_2=\frac{1}{|Q|}\int_Q|T_Df_2|dx<\left(\frac{1}{|Q|}\int_Q|T_Df_2|^2dx\right)^{1/2}.$$
By boundedness of a smooth convolution Calder\'{o}n-Zygmund operator in $L^2$ and by triangle inequality we obtain
$$I_2\lesssim_{T,D}\left(\frac{1}{|Q|}\int_{2Q}|f_2|^2dx\right)^{1/2}=
\left(\frac{1}{|Q|}\int_{2Q}|f-f_Q|^2dx\right)^{1/2}$$
$$\leq\left(\frac{1}{|Q|}\int_{2Q}|f-f_{2Q}|^2dx\right)^{1/2}
+\left(\frac{1}{|Q|}\int_{2Q}|f_Q-f_{2Q}|^2dx\right)^{1/2}.$$
The first summand is estimated by inverse H\"{o}lder inequality
$$\lesssim\frac{1}{|Q|}\int_{2Q}|f-f_{2Q}|dx\lesssim \omega(\ell)$$
(see Introduction). For the second summand we use   (\ref{eq:eq9}) and we are done for $k=2.$

b) Taking care of third term is not so easy. To estimate oscillation
$$I_3=\frac{1}{|Q|}\int_Q|Tf_3-b_{3,Q}|dx$$
we choose an appropriate constant $b_{3,Q}=Tf_3(x_0),$ where $x_0$ is a centre of $Q.$
A smooth  kernel $K$ satisfies an estimate
$$\nabla K(x_0)\lesssim \frac{1}{|x_0|^{d+1}}, \quad x_0\neq 0.$$
This implies
$$|K(u-x)-K(u-x_0)|$$
$$\lesssim \sup _{|t-x_0|<|x-x_0|}|\nabla K(u-t)||x-x_0|$$

$$\lesssim \frac{|x-x_0|}{|u-t|^{d+1}},$$
with the constants not depending on $u, x, x_0$ and $Q.$
Substituting this to estimate $I_3$ we have
$$I_3\lesssim \frac{C}{|Q|}\int_Q dx \int_{D\backslash 2Q} \frac{|x-x_0|}{|u-t|^{d+1}}|f-f_Q|(u)du$$
$$ \lesssim \ell\int_{D\backslash 2Q} \frac{|f-f_Q|(u)}{|u-t|^{d+1}}du.$$
Now, replace $f$ by its extension $\tilde{f}$  and note that $f_Q=\tilde{f}_Q$ for $Q\subset D$ we produce:
 $$ I_3\lesssim \ell\int_{\mathbb{R}^d\backslash 2Q} \frac{|\tilde{f}-\tilde{f}_Q|(u)}{|u-t|^{d+1}}du.$$
Put $Q_k=2^kQ\backslash2^{k-1}Q,$ and  split integral in sum as following
$$I_3\lesssim\ell\sum_{k=1}^\infty \frac{1}{(\ell2^k)^{d+1}}\int_{Q_k} |\tilde{f}-\tilde{f}_Q|(u)du.$$
Applying the telescopic summation method we have
$$I_3\lesssim\ell\sum_{k=1}^\infty \frac{1}{(\ell2^k)^{d+1}}\sum_{s=0}^k \int_{Q_k} |\tilde{f}_{2^sQ}-\tilde{f}_{2^{s+1}Q}|(u)du,$$
where we denoted $\tilde{f}_{2^{k+1}Q}=\tilde{f}$ in each inner sum.
 Continue calculations with the help of (\ref{eq:eq9}), we obtain
 $$I_3\lesssim\ell\sum_{k=1}^\infty \frac{1}{(\ell2^k)^{d+1}}\sum_{s=0}^k |Q_k| |\tilde{f}_{2^sQ}-\tilde{f}_{2^{s+1}Q}|\lesssim$$
 $$\lesssim\ell\sum_{k=1}^\infty \frac{1}{(\ell2^k)^{d+1}}\sum_{s=0}^k (\ell 2^k)^d \omega (\ell 2^s)\|f\|_{\omega,D}\approx$$
$$\approx\sum_{k=1}^\infty \frac{1}{2^k}\sum_{s=0}^k  \omega(2^s\ell)\|f\|_{\omega,D}.$$
Changing the order of summing continue
$$I_3\lesssim \sum_{s=1}^\infty\frac{1}{2^s}\omega(2^s\ell)\|f\|_{\omega,D}$$
$$\lesssim \int_1^\infty \frac{\omega(t\ell)}{t^2}dt.$$
Put a change of variable, take into account that  the function $\frac{\omega(t)}{t^{\varepsilon}}$ is almost decreasing for some $\varepsilon<1$ we have
$$I_3\lesssim_{T,D,f}\omega(\ell),$$
and we are done.$\Box$
\end{pf}

In fact, we proved that  (\ref{eq:eq8})  is equivalent to inequality
\begin{equation}
\label{eq:eq10}
  \frac{1}{|Q|}\int_{Q }|T_D f_1-b_Q|dx=\frac{|f_Q|}{|Q|}\int_{Q }|T \chi_D-b_Q|dx\lesssim \omega(\ell)\|f\|
\end{equation}
Now we will estimate mean values  $f_Q$ over cubes $Q.$
\begin{lem}
\label{lem5}
   There exists a  constant $C>0,$ depending on $D,$ such that for any  $f\in\mathcal{C}_{\omega}(D)$ and each cube $Q\subset D$ one has

   $$|f_Q|\leq C \|f\|\int_{\ell}^1\frac{\omega(t)}{t}dt.$$

   \end{lem}
   \begin{pf}
Choose minimal integer $n$ such that $2^n\ell>R, $  replace $f$ by its extension $\tilde{f}$  and note that $f_Q=\tilde{f}_Q$ for $Q\subset D.$ Then by telescopic summation by means of (\ref{eq:eq9}) one produces
$$|f_Q|\leq\sum_0^{n-1}|\tilde{f}_{2^sQ}-\tilde{f}_{2^{s+1}Q}|+|\tilde{f}_{2^n Q}|\lesssim$$
$$\lesssim\sum_0^{n-1}\omega(2^s\ell)\|\tilde{f}\|_{\omega,D}+\|\tilde{f}\|_{L^1(\mathbb{R}^d)}\lesssim$$
$$\lesssim \|\tilde{f}\|\int_{\ell}^1\frac{\omega(t)}{t}dt$$
$$ \lesssim \|f\|\int_{\ell}^1\frac{\omega(t)}{t}dt ,$$
that is claimed. $\Box$
\end{pf}

If $\omega$ is Dini regular this implies that all mean values  $|f_Q|$ are uniformly bounded.
In this case we obtain exactly T1 theorem. If $\omega$ fails to be  Dini regular, then $(T\chi_D)\chi_D \in \mathcal{C}_{\tilde{\omega}}(D),$ is sufficient for the  boundedness of the restricted $T_D$ in the space $\mathcal{C}_{\omega}(D).$ Therefore the sufficienty condition of Theorem \ref{thm2} is proved.

In order to prove the necessity condition, we  construct a  function with extremal properties. Define a radial function in $\mathbb{R}^d$ by formula $\varphi(t)=\int_{|t|}^1\frac{\omega(u)}{u}du$ when $|t|<1.$ We extend this function outside the open ball $|t|<1$ by zero. We need a modification of Lemma 1.1 \cite{Sj}.
\begin{lem}
\label{lem6}
   We have $\varphi\in\mathcal{C}_{\omega}(\mathbb{R}^d).$
   For arbitrary  cube $Q$ with the centre in $0$  and side length $\ell <1,$   mean values of $\varphi$
   are bounded from below
   $$|\varphi_Q|\gtrsim\int_{\ell}^1\frac{\omega(t)}{t}dt$$
with a constant not depending on $Q.$
   \end{lem}
\begin{pf}
It is easy exercise that we can use integration over balls instead of cubes in
definition of Campanato spaces.  So, we may calculate a mean value $\varphi_B$  for the ball $B$ with a centre $0$ and radius $\ell,$ and obtain that  $\varphi_B\approx\int_{\ell}^1\frac{\omega(u)}{u}du.$
The same estimate takes place for cubes and we are done. $\Box$
\end{pf}

Now we are ready to finish the proof of necessity. Given $\tau\in D$ define the function  $\varphi_\tau (t)=\varphi (t-\tau).$ In the Lipschitz domain $D$  the function  $\varphi_\tau\chi_D$ has the properties similar of  $\varphi.$ That is
   $\varphi_\tau\chi_D\in\mathcal{C}_{\omega}(D)$ with norms in $\mathcal{C}_{\omega}(D)$ uniformly bounded in $\tau\in D$ by a constant depending on the Lipschitz constant of the boundary of $D.$   Also,  the mean values satisfy
   \begin{equation}\label{eq:eq11}
  (\varphi_\tau\chi_D)_Q\gtrsim\int_{\ell}^1\frac{\omega(t)}{t}dt
\end{equation}
with constant not depending on $Q$ and $\tau.$
Finally, we substutute  $f=\varphi_\tau$ in (\ref{eq:eq10}), then  (\ref{eq:eq11}) implies there exist a constant $b_Q,$ such that
$$\frac{1}{|Q|}\int_{Q }|T \chi_D-b_Q|dx\lesssim \tilde{\omega}(\ell),$$
The proof of Theorem \ref{thm2} is finished. $\Box$

\section {Proof of  Theorem 3.}
 We start with imbedding of the Bloch space into the Campanato space.
\begin{lem}
\label{lem7}
 Let $\omega$ be a   modulus of continuity   and let $D\subset \mathbb{R}^d$ be a  bounded Lipschitz domain. Then $ \mathcal{B}_{\omega} (D)\subset\mathcal{C}_{\omega} (D).$
   \end{lem}
\begin{pf}
 An Extension Lemma \ref{lem1}  implies that we may consider only cubes $Q$ such that $2Q\subset D.$
We take care in  estimate for arbitrary $f \in \mathcal{B}_{\omega} (D)$ and appropriate constant $b_Q$

$$I=\frac{1}{|Q|}\int_Q|f-b_Q|dx\leq C  \omega(\ell),$$
 where a constant $C$ does not depend on $f$ and $Q.$ Given  a centre of the cube  $x_0 \in D$  we take a constant $b_Q=f(x_0).$ Then easy calculations give us
 $$f(x)-f(x_0)=\int_0^1\langle \nabla f(x_0+t(x-x_0),x-x_0\rangle dt.$$

 Take spherical coordinates  a pole in $x_0$ one obtains
$$\frac{1}{|Q|}\int_Q|f(x)-f(x_0)|dx\lesssim\frac{1}{|Q|}\int_{S^d} dS^d \int_{o}^{2\ell}r^{d-1}dr\int_0^1|\langle \nabla f(x_0+t(x-x_0),x-x_0\rangle |dt$$
$$\lesssim\frac{1}{|Q|}\int_{o}^{2\ell}r^{d-1}dr\omega(\ell)\lesssim \omega(\ell),$$
 that is claimed.
 \end{pf}
  Now we prove that this imbedding is sharp. Define the subspaces $\mathcal{C}^h_{\omega} (D)\in \mathcal{C}_{\omega} (D)$  and  $\mathcal{B}^h_{\omega} (D)\in \mathcal{B}_{\omega} (D,)$ consisting of harmonic functions in $D.$

 \begin{lem}
\label{lem8}

   Let $\omega$ be a   modulus of continuity   and let $D\subset \mathbb{R}^d$ be a  bounded Lipschitz domain domain.  Then $ \mathcal{C}^h_{\omega} (D)\subset\mathcal{B}^h_{\omega} (D).$
   \end{lem}
\begin{pf}
 It is a simple modification of a  result for analytic Bloch class  in the disk.
 Fix a ball $B=B(x_0,R)$ with a centre in $x_0 \in D$ and radius $R< \frac{1}{2}dist(x, \partial D).$
 Represent a function $f\in \mathcal{C}^h_{\omega} (D)$ in a point $x\in B$ by Poisson formula
 $$f(x)= \frac{C_d}{r}\int_{|y-x|=r} \frac{r^2-|x-x_0|^2}{|y-x|^d}f(y)dS_r(y)$$
 with $|x-x_0|<r< R$  and $dS$ is induced surface masure. Now for arbitrary constant $C$ we have
 $$f(x)-C= \frac{C_d}{r}\int_{|y-x|=r} \frac{r^2-|x-x_0|^2}{|y-x|^d}(f(y)-C)dS_r(y)$$
Then differentiate in the point $x_0$ we reduce
$$\frac{\partial}{\partial x_i}f(x_0)= \frac{C_d}{r}\int_{|y-x|=r} \frac{\partial}{\partial x_i}\left(\frac{r^2-|x-x_0|^2}{|y-x|^d}\right)|_{x=x_0}(f(y)-C)dS_r(y).$$
Easy estimates of the derivatives of the Poisson kernel show that
$$\frac{\partial}{\partial x_i}f(x)\lesssim  \frac{1}{r^{d}}\int_{|y-x|=r} |f(y)-C|dS_r(y).$$
The left part does not depend on $r$ and  after integration on $r$ one obtains

$$\frac{\partial}{\partial x_i}f(x_0)\lesssim \frac{1}{R}\int_R^{2R} \frac{dr}{r^{d}}\int_{|y-x|=r} |f(y)-C|dS_r(y)$$
$$\leq\frac{1}{R^{d+1}}\int_0^{2R} dr\int_{|y-x|=r} |f(y)-C|dS_r(y)$$
$$\leq\frac{1}{R^{d+1}}\int_{B(x_0,2R)}  |f(y)-C|dy.$$
Now we may take $C$ being a mean value over the ball $B(x_0,2R)$ and  finish the proof by the same argument as in Lemma \ref{lem6}.
 \end{pf}

 From now on we consider smooth domains, so that the number of properties of the following definition is satisfied.
 \begin{df}
 \label{df4}
 Given a modulus of continuity $\omega,$ a bounded domain $D\subset  \mathbb{R}^d$ is called to be a $C^{1,\omega}-$smooth domain if there is exist positive constant $r_0,$ such that for  each  $a\in\partial D$ there exist the function $A \in C^{1,\omega}(\mathbb{R}^{d-1})$  and a cube $Q$ with side length $r_0$ and centre $a$ such that, after a suitable shift and rotation
 $$D\cap Q=\{ (x',x_d)\in (\mathbb{R}^{d-1}, \mathbb{R})\cap Q:x_d>A(x')\}.$$
Under this rotation we may assume that the tangent hyperplane in $a=0$ to $\partial\Omega$ is $x_d=0,$  and $\nabla A(0')=0,$ Moreover,

  \begin{equation}\label{eq:eq12}
    |A(x')|\leq c_0 |x'|\omega(|x'|), |x'|<r_0 .
  \end{equation}
\end{df}
We start with a formula of differentiation of $PV$ integrals.
  \begin{lem}
\label{lem9}
  Recall that $ K(x)$
 being a homogeneous function of degree $-d$, is continuously differentiable  on
$\mathbb{R}^d \setminus 0.$  Then for $y$ in $C^{1,\omega}-$smooth domain $D:$

\begin{equation}\label{eq:eq13}
\frac{\partial}{\partial y_i} PV\int_D K(x-y)dx=-\int_{\partial D}K(x-y) \cos(\nu,x_i)dS(x),
\end{equation}
provided $\nu$ is  outer normal to $\partial D,$ and $dS(x)$ is an induced  surface measure on $\partial D.$
\end{lem}
\begin{pf}
Consider a formula (3) \S8, Ch.2 in \cite{Mikh} of differentiation under integral sign,
\begin{equation}
\label{eq:eq14}
\frac{\partial}{\partial y_i} \int_{D_\varepsilon} K(x-y)dx=  \int_{D_\varepsilon} \frac{\partial}{\partial y_i} K(x-y)dx +\int_{S_\varepsilon}K(x-y) \cos(\nu,x_i)dS(x),
\end{equation}
where we denote $D_\varepsilon=D\setminus \{x: |x-y|<\varepsilon\},\:S_\varepsilon=\{x: |x-y|=\varepsilon\}$ and $\varepsilon< \textrm{dist}(y,\partial D).$ By Green formula applied to the first integral on the right part of (\ref{eq:eq14}) one has
 $$\int_{D_\varepsilon} \frac{\partial}{\partial y_i} K(x-y)dx= -\int_{D_\varepsilon} \frac{\partial}{\partial x_i} K(x-y)dx$$
 $$= -\int_{\partial D \cup S_\varepsilon}K(x-y) \cos(\nu,x_i)dS(x).$$
So, we obtain
$$\frac{\partial}{\partial y_i} \int_{D_\varepsilon} K(x-y)dx=-\int_{\partial D}K(x-y) \cos(\nu,x_i)dS(x).$$
The right part does not depend on $\varepsilon,$ and  by theorem on differentiation with respect to parameter we are done  for any $y\in D.$
\end{pf}

There is a
lemma, which plays a crucial role for the study of the  smooth convolution Calder\'{o}n-Zygmund operator with even kernel  \cite{I,MOV}.
\begin{lem}
\label{lem10}(Extra cancellation property)
  Let   $B$ be an arbitrary ball in $\mathbb{R}^d.$ Then $(T\chi_B)|_B\equiv 0.$

\end{lem}
As a corollary, we have
\begin{equation}
 \label{}
    (\frac{\partial}{\partial x_i}  T\chi_B)\chi_B\equiv 0.
 \end{equation}

 To check the condition $(T\chi_D)\chi_D\in \mathcal{B}_{\omega} (D),$ we use a modification of Theorem 1 in \cite{V}.
Let  $y\in D.$ Put  $\delta=\textrm{dist}(y,\partial lD).$ Let $a\in\partial D$ be the point where this distance is attained. Applying Definition \ref{df4}, we may assume  that $a=0,$ and   $(0',t)\in D$ when $0<t<\delta$ and  $y=(0',\delta)$

Let $B=B((0',r_0), r_0)$ be the ball with a centre $(0',r_0)$ and radius $r_0$ tangent to $\partial D$ in $0$. By Lemma \ref{lem10}, we  proceed
$$\frac{\partial}{\partial y_i} T\chi_D= \frac{\partial}{\partial y_i} (T\chi_D- T\chi_B)$$
$$=\frac{\partial}{\partial y_i} (T\chi_{D\backslash B}-  T\chi_{B\backslash D}).$$
And now by Lemma \ref{lem9} and Green formula we continue
$$\frac{\partial}{\partial y_i} T\chi_D=\int_{B\backslash D} \frac{\partial}{\partial y_i} K(x-y)dx -\int_{B\backslash D} \frac{\partial}{\partial y_i} K(x-y)dx.$$

 The obvious advantage is that domain is tangent to $\partial D$ in $a=0$ and hence is small. Further we estimate the absolute value of integral.
 $$I=\int_{B\backslash D \cup B\backslash D} \left |\frac{\partial}{\partial y_i} K(x-y)\right|dx .$$
 We split $I$ over the level $r_0$
 $$I=\int_{|x|<r_0} + \int_{|x|>r_0} .$$
 Clearly, for the second integral one has
 $$ \int_{|x|>r_0} =\int_{(B\backslash D \cup B\backslash D)\cap\{|x|>r_0\}} \left|\frac{\partial}{\partial y_i} K(x-y)\right|dx $$
 $$ \lesssim \frac{1}{r_0^{d+1}} |D|, $$
 which does not depend on $\delta.$

  We  split the first integral over level $ |x'|=\delta$ and proceed
   $$ \int_{|x|<r_0} =\int_{|x|<\delta}+\int_{\delta< |x|< r_0}=II+III$$
 In order to estimate  term $II$ note that $|x|<\delta$ implies that $|x'|^2+(x_d-\delta)^2\gtrsim \delta^2.$
 Hence, we have
 $$II\lesssim\int_{|x'|<\delta}\frac{dx'}{\delta^{d+2}}\int_{|x_d|\lesssim |x'|\omega(|x'|)}dx_d$$
 $$\lesssim\int_{|x'|<\delta}\frac{|x'|\omega(|x'|)dx'}{\delta^{d+2}}$$
 $$\lesssim \frac{\omega(\delta)}{\delta}.$$

 For the term $III$ note that $\delta< |x|<r_0$ implies that $|x'|^2+(x_d-\delta)^2\gtrsim |x'|^2.$
 Hence
$$III\lesssim \int_{\delta< |x'|< r_0}\frac{dx'}{|x'|^{d+2}}\int_{|x_d|\lesssim |x'|\omega(|x'|)}dx_d$$
 $$\lesssim\int_{\delta< |x'|< r_0}\frac{\omega(|x'|)dx'}{|x'|^{d+1}}$$
$$\lesssim \frac{\omega(\delta)}{\delta},$$
because  $\frac{\omega(t)}{t^{\varepsilon}}$ is almost decreasing for some $\varepsilon<1.$
 This completes the proof of Proposition  \ref{pr2}.$\Box$

 Therefore, the last reduction of Theorem \ref{thm3} is proved.




\end{document}